\documentclass{mathincs}

\usepackage{hyperref}
\usepackage{pict2e}
\usepackage{ketpic2e}
\usepackage{amsmath}
\usepackage{amsthm}
\usepackage{graphicx}
\usepackage{animate}
\usepackage{xcolor}
\definecolor{labelcolor}{RGB}{100,0,0}
 \theoremstyle{definition}
 
 \theoremstyle{remark}

 \numberwithin{equation}{section}

\makeatletter
\newcommand{\verbatimfont}[1]{\renewcommand{\verbatim@font}{\ttfamily#1}}
\makeatother

\begin{document}
\verbatimfont{\small}
%
%
%
\title[Hamiltonian systems: symbolical, numerical and 
graphical study]{Hamiltonian dynamical systems: symbolical, numerical and 
graphical study}
\author[Takato]{Setsuo Takato}

\address{Faculty of Sciences\\
T\={o}h\={o} University\\
Funabasi (Chiba) Japan}

\email{takato@phar.toho-u.ac.jp}

\thanks{Work partially supported by grants KAKENHI 15K01037 (ST, Japan) and CONACyT CB-2012 179115 (JAV, M\'exico).}
\author[Vallejo]{Jos\'e A Vallejo}
\address{Faculty of Sciencies\br
State University of San Luis Potos\'i\br
San Luis Potos\'i (SLP) M\'exico}
\email{jvallejo@fc.uaslp.mx}


\date{\today}

\begin{abstract}
Hamiltonian dynamical systems can be studied from a variety of viewpoints.
Our intention in this paper is to show some examples of usage of two Maxima
packages for symbolical and numerical analysis (\texttt{pdynamics} and
\texttt{poincare}, respectively), along with the set of
scripts \KeTCindy\ for obtaining the \LaTeX\ code corresponding to graphical representations
of Poincar\'e sections, including animation movies.  
\end{abstract}

\maketitle
\section{Introduction}

For simplicity, we will consider Hamiltonians defined on the symplectic manifold
$\mathbb{R}^{2n}$, with coordinates $(q^j,p_j)$ ($1\leq j\leq n$), endowed with 
the canonical form $w=\mathrm{d}p_j\wedge \mathrm{d}q^j$, and the induced Poisson bracket
on $\mathcal{C}^\infty (\mathbb{R}^{2n})$
$$
\{f,g\}=\sum^n_{i=1}\left( \frac{\partial f}{\partial p_i}\frac{\partial f}{\partial q^i}
-\frac{\partial f}{\partial q^i}\frac{\partial f}{\partial p_i}\right)\,,
$$
although all the results remain valid for an arbitrary symplectic manifold. 
For background on Hamiltonian systems, see \cite{Cus97}.

Given a Hamiltonian system defined by the Hamiltonian function 
$H\in\mathcal{C}^\infty (\mathbb{R}^{2n})$ and the Poisson bracket $\{\cdot ,\cdot\}$,
\begin{align}\label{eq1}
\dot{q}^j &= \frac{\partial H}{\partial p_j} \nonumber\\
\dot{p}_j &= -\frac{\partial H}{\partial q^j}\,,
\end{align}
two of the main goals in the theory of dynamical systems are the determination of
possible closed, stable orbits, and the computation of adiabatic invariants (of
course, taking for granted the impossibility of solving \eqref{eq1} explicitly).
Of particular interest is the case in which the Hamiltonian $H$ is a perturbation of
an integrable one, say, $H=H_0 +\sum^n_{j=1}\varepsilon^j H_j$. A widely used
procedure to study it, consists in writing the Hamiltonian in the so-called
\emph{normal form}, that is, as a formal series
\begin{equation*}
H=\sum^\infty_{j=0}\varepsilon^j N_j
\end{equation*}
where $N_0=H_0$, and each $N_j$ commutes with the unperturbed Hamiltonian,
$$
\{H_0,N_j\}=0\,.
$$

Let us recall that given an integral curve, that is, a $c(t)=(q(t),p(t))$ satisfying \eqref{eq1},
the evolution of any observable $f\in \mathcal{C}^\infty (\mathbb{R}^{2n})$ along $c$ is determined by
$$
\dot{f}(t)=\{H, f\}\,.
$$
Any smooth function such that $\{H,f\}=0$ is thus a constant of motion, also called a first integral.
Indeed, given enough first integrals one can solve the motion of the system, as the physical
trajectories are determined by the intersection of their level hypersurfaces. Unfortunately, determining
first integrals is a very difficult problem, and there are quite a few systems for which enough
first integrals exist (roughly, these are the so-called integrable systems).

Notice that transforming to the normal form introduces a (possibly infinite)
family of first integrals $N_j$ which might not been present in the original system. 
These additional, spurious symmetries must be removed in order to have a system equivalent to
the original one, and this is usually done by
restricting the system to a reduced phase space through symplectic (singular)
reduction. The basic idea is to restrict the system to a particular level hypersurface and
to consider its evolution there.
A number of well-known techniques are available to do this, for instance the
ones based on Moser's theorem \cite{Mos70}: If $M_h$ denotes the hypersurface
$H_0=h$, suppose that the orbits of the Hamiltonian flow $\mathrm{Fl}^t_{X_{H_0}}$ are all periodic 
with period $T$ and let $S$ be the quotient with respect to the
induced $U(1)-$action on $M_h$. Then, to every non-degenerate critical point
$\overline{p}\in S$ of the restricted averaged perturbation
$\left. N_1\right|_S=\left. \left\langle H_1\right\rangle\right|_S$
corresponds a \emph{periodic} trajectory of the full Hamiltonian vector field $X_H$, that branches
off from the orbit represented by $\overline{p}$ and has period close to $2\pi$.
When the critical points are degenerate, one can resort to the second-order normal
form to decide the stability of orbits. An example of this situation is given by the
H\'enon-Heiles Hamiltonian \cite{Cus94}. These results illustrate the importance of
being able to compute efficiently the normal form of a Hamiltonian system. In Section
\ref{sec2} we show how to do this using the Computer Algebra System (CAS) Maxima.

Another aspect related to the study of existence and stability of closed orbits is
the construction of Poincar\'e sections. They provide a direct and very intuitive
way for detecting these orbits, but their computation in closed
form is usually impossible, so numerical methods are needed. The traditional method
used for this task has been the fourth-order Runge-Kutta, but more recently methods
based on symplectic integrators (such as symplectic Euler, St\"ormer-Verlet, symplectic
RK, etc.) are also intensively used, see \cite{BC16} for a recent review.
The choice of one method or another depends very
much on the properties of the system under consideration, in Section \ref{sec3} we will use the RK method, but symplectic methods can be included by substituting the \texttt{rkfun} command in the code of
\texttt{poincare} with \texttt{symplectic\_ode}, recently included in Maxima
(from version 5.39.1 onward). In any case, one of the main
goals is to obtain a clean picture of the phase-space portrait of the system, something
that can be challenging for CASs, whose graphical
output is not very sophisticate in many cases. To deal with this issue we present in
Section \ref{sec4} the set of CindyScript macros called \KeTCindy\,, 
which parse the output of Maxima through
the Dynamical Geometry Software (DGS) Cinderella and return the \LaTeX\  code of the
corresponding graphics, that can be included in any document even in the form of an animation.
The data for these graphics are actually codes of TPIC specials for \LaTeX\,, or
pict2e commands in the case of pdf\LaTeX\,, so there they can be inserted in scientific
documentation with great flexibility (see \cite{icms2016,cas2016,iccsa2017,castr2017}
for installation and examples of use).

The Maxima packages \texttt{pdynamics} (short for `Poisson Dynamics') and
\texttt{poincare} can be downloaded from \url{http://galia.fc.uaslp.mx/~jvallejo/pdynamics.zip}
and \url{http://galia.fc.uaslp.mx/~jvallejo/poincare.mac} respectively\footnote{There is a
documentation file inside \texttt{pdynamics.zip}, and the documentation for \texttt{poincare.mac}
can be found at \url{http://galia.fc.uaslp.mx/~jvallejo/PoincareDocumentation.pdf}. Both files
contain detailed instructions about the installation.}. The
\KeTCindy\ package is available at \url{http://ketpic.com/?lang=english}.

\section{Symbolic study of Hamiltonian systems: normal forms}\label{sec2}
Given a smooth vector field in $\mathbb{R}^{m}$ (although what follows is valid in an
arbitrary manifold) its flow is a mapping $\mathrm{Fl}_X:\mathbb{R}^{m}\to\mathbb{R}^{m}$
defined by
$$
\mathrm{Fl}_X(t,p)\doteq\mathrm{Fl}^t_X(p)\doteq c_p(t)\,,
$$
where $c_p$ is the integral curve of $X$ such that, for $t=0$, passes through $p\in\mathrm{R}^m$
(i.e., $c_p(0)=p$). When $m=2n$, there is a canonical symplectic form
$$
\Omega = \mathrm{d}p_1\wedge\mathrm{d}q_1+\cdots +\mathrm{d}p_n\wedge\mathrm{d}q_n\,.
$$

Any Hamiltonian $H\in\mathcal{C}^\infty (\mathbb{R}^{2n})$, has an associated vector field
$X_H$ defined by the condition $i_{X_H}\Omega =-\mathrm{d}H$. In local coordinates
$(q_i,p_i)$ it has the expression
$$
X_H =\left( \frac{\partial H}{\partial p_1},-\frac{\partial H}{\partial q_1},\ldots,
\frac{\partial H}{\partial p_n},-\frac{\partial H}{\partial q_n}\right)\,.
$$ 

Suppose now that $X$ is the generator of an $\mathbb{S}^1-$action, so the flow $\mathrm{Fl}^t_X$
is periodic in the variable $t$. This property can be used to put $H$ in normal form (for details,
see \cite{AVV13}). To this end, it is essential to define two averaging operators acting on
observables. The first one is denoted as $\left\langle\cdot\right\rangle$ and is given by
integrating the pullback
$$
\left\langle g\right\rangle\doteq \frac{1}{2\pi}\int^{2\pi}_0(\mathrm{Fl}^t_X)^*g\,\mathrm{d}t\,,
$$
for any observable $g\in\mathcal{C}^\infty (\mathbb{R}^{2n})$. The second operator, denoted
$\mathcal{S}$, is defined as
$$
\mathcal{S}(g)\doteq\frac{1}{2\pi}\int^{2\pi}_0(t-\pi)(\mathrm{Fl}^t_X)^*g\,\mathrm{d}t\,. 
$$
In the particular case of a perturbed Hamiltonian, of the form $H=H_0 +\varepsilon H_1
+\frac{\varepsilon^2}{2}H_2 +\cdots$, if the non-perturbed part $H_0$ generates an $\mathbb{S}^1-$action in such a way that its flow is periodic with frequency function
$w$, it can be proved (see \cite{AVV13}) that its second-order normal form is
$$
N=H_0 +\varepsilon\left\langle H_1\right\rangle +
\frac{\varepsilon^2}{2}\left(\left\langle H_2\right\rangle +
\left\langle\left\lbrace\mathcal{S}\left(\frac{H_1}{w}\right) ,H_1\right\rbrace\right\rangle
 \right)\,.
$$

There are other representations of the normal form (it must be stressed that it is \emph{not}
unique), but this one has the particular features that it is global (not depending on action-angle
variables), and particularly well-suited for symbolic computation. Let us illustrate the use of
the \texttt{pdynamics} Maxima package by considering the example of the Pais-Uhlenbeck oscillator.
This system is a toy model of a field theory defined by a Lagrangian depending on higher-order 
derivatives. These Lagrangians are believed to lead to perturbatively renormalizable
theories, where the infinities appearing in the perturbation series for the field equations can
be cured through some well-defined regularization procedure. The corresponding Hamiltonian is constructed through a higher-order analog of the Legendre transformation, called the
Ostrogadskii formalism \cite{Os50}. After some suitable transformations (see \cite{Pav13})
the Hamiltonian can be expressed
as the \emph{difference} of two harmonic oscillators with respective frequencies $w_1$, $w_2$.
Adding an interaction term in the form of a homogeneous polynomial results in the Hamiltonian
\begin{equation}\label{puham}
H=\frac{1}{2}(p^2_1 +w^2_1q^2_1)-\frac{1}{2}(p^2_2 +w^2_2q^2_2)+\frac{\lambda}{4}(q_1+q_2)^4\,.
\end{equation}
This can be considered as a perturbed system of the form $H=H_0+\lambda H_1$,
let us study it symbolically. We would use the following sequence of commands in Maxima:

\noindent
\begin{minipage}[t]{8ex}\bf
(\% i1) 
\end{minipage}
\begin{minipage}[t]{\textwidth}\begin{verbatim}
load(pdynamics)$\end{verbatim}
\end{minipage}

\noindent
\begin{minipage}[t]{8ex}\bf
(\% i2)
\end{minipage}
\begin{minipage}[t]{\textwidth}\begin{verbatim}
declare(w1,integer)$\end{verbatim}
\end{minipage}

\noindent
\begin{minipage}[t]{8ex}\bf
(\% i3)
\end{minipage}
\begin{minipage}[t]{\textwidth}\begin{verbatim}
assume(w1>0)$\end{verbatim}
\end{minipage}

\noindent
\begin{minipage}[t]{8ex}\bf
(\% i4)
\end{minipage}
\begin{minipage}[t]{\textwidth}\begin{verbatim}
declare(w2,integer)$\end{verbatim}
\end{minipage}

\noindent
\begin{minipage}[t]{8ex}\bf
(\% i5)
\end{minipage}
\begin{minipage}[t]{\textwidth}\begin{verbatim}
assume(w2>0)$\end{verbatim}
\end{minipage}

\noindent
\begin{minipage}[t]{8ex}\bf
(\% i6)
\end{minipage}
\begin{minipage}[t]{\textwidth}\begin{verbatim}
H0(q1,p1,q2,p2):=(p1^2+w1^2*q1^2)/2-(p2^2+w2^2*q2^2)/2$\end{verbatim}
\end{minipage}

\noindent
\begin{minipage}[t]{8ex}\bf
(\% i7)
\end{minipage}
\begin{minipage}[t]{\textwidth}\begin{verbatim}
H1(q1,p1,q2,p2):=(q1+q2)^4/4$\end{verbatim}
\end{minipage}

\noindent
\begin{minipage}[t]{8ex}\bf
(\% i8)
\end{minipage}
\begin{minipage}[t]{\textwidth}\begin{verbatim}
H2(q1,p1,q2,p2):=0$\end{verbatim}
\end{minipage}\\

Up to here, we have just defined the parameters of the system (the frequencies $w_1$, $w_2$, and
the subhamiltonians $H_i$). Let us check that the Hamiltonian flow of the non-perturbed part $H_0$ is
periodic by explicitly computing it (we have slightly edited the output of Maxima by writing it as a
column matrix, to make it more readable):

\noindent
\begin{minipage}[t]{8ex}\bf
(\% i9)
\end{minipage}
\begin{minipage}[t]{\textwidth}\begin{verbatim}
phamflow(H0);\end{verbatim}
\end{minipage}
\[\tag{\% o9}
\begin{pmatrix}
\frac{p_1\, \sin{\left( t\, w_1\right) }}{w_1}+q_1\, \cos{\left( t\, w_1\right) }\\
p_1\, \cos{\left( t\, w_1\right) }-q_1\, w_1\, \sin{\left( t\, w_1\right) }\\
q_2\, \cos{\left( t\, w_2\right) }-\frac{p_2\, \sin{\left( t\,w_2\right) }}{w_2}\\
q_2\, w_2\, \sin{\left( t\, w_2\right) }+p_2\, \cos{\left( t\, w_2\right) }
\end{pmatrix}
\] \\

It is clear that the flow is periodic with period $T=2\pi w_1w_2$. Thus, we define the
frequency function as

\noindent
\begin{minipage}[t]{8ex}\bf
(\% i10)
\end{minipage}
\begin{minipage}[t]{\textwidth}\begin{verbatim}
u(q1,p1,q2,p2):=1/(w1*w2)$\end{verbatim}
\end{minipage}\\

Finally, we can compute $N_1=\left\langle H_1\right\rangle$:

\noindent
\begin{minipage}[t]{8ex}\bf
(\% i11)
\end{minipage}
\begin{minipage}[t]{\textwidth}\begin{verbatim}
phamaverage(H1,H0,u(q1,p1,q2,p2));\end{verbatim}
\end{minipage}
\begin{align*}
\tag{\% o11} 
\frac{1}{32 {w_1^{4}}\, {w_2^{4}}}\left[\left( \left( 3 {q_2^{4}}+12 {q_1^{2}}\, {q_2^{2}}+3 {q_1^{4}}\right) \, {w_1^{4}}+\left( 12 {p_1^{2}}\, {q_2^{2}}+6 {p_1^{2}}\, {q_1^{2}}\right) \, {w_1^{2}}+3 {p_1^{4}}\right) \, {w_2^{4}}\right.\\
\left. +\left( \left( 6 {p_2^{2}}\, {q_2^{2}}+12 {p_2^{2}}\, {q_1^{2}}\right) \, {w_1^{4}}+12 {p_1^{2}}\, {p_2^{2}}\, {w_1^{2}}\right) \, {w_2^{2}}+3 {p_2^{4}}\, {w_1^{4}}\right]
\end{align*}

The second-order normal form can be computed along similar lines but, as one can guess, the
expressions become very cumbersome, and not much illuminating
(see \eqref{o15} below). Indeed, it is customary to
simplify these expressions by rewriting them in terms of the so-called Hopf variables. The
idea behind these variables is the following: in the case in which the normal subhamiltonians
$N_i$ are polynomials, the fact that they commute with $N_0=H_0$ means that they are invariant under the smooth $\mathcal{S}^1-$action of $X_{H_0}$. The space of smooth invariant functions
is finitely generated (this is a generalization to the smooth case of a classical result of
Hilbert dealing with algebraic invariants, called the Schwarz theorem \cite{Sch75}), and a set of functional generators is precisely given by the Hopf polynomials, that can be considered as new variables.
In other words, any smooth invariant function can be expressed as a smooth function of the Hopf variables. For the Pais-Uhlenbeck oscillator with resonance $1:2$ (that is, when $w_1=1$ and
$w_2=2$), these Hopf invariants can be readily computed \cite{AVV17} and they turn out to be
\begin{align}\label{rhos}
\rho_1 =& q^2_1+p^2_1\nonumber \\
\rho_2=& 4 q^2_2+p^2_2\nonumber \\
\rho_3 =& p_2(p^2_1-q^2_1)-4p_1q_1q_2 \\
\rho_4 =& 2q_2(p^2_1-q^2_1)+2q_1p_1p_2\nonumber \,.
\end{align}

We can compute $N_2$ by extracting the coefficient of $\lambda^2$ in the second-order normal form
of $H$. The following commands show how to study this resonance, defining a function
\texttt{phopf6res12} (not contained in the \texttt{pdynamics} package) adapted to this case, whose
purpose is to express everything in terms of the variables \eqref{rhos}. First, we define the
Hamiltonian:

\noindent
\begin{minipage}[t]{8ex}\bf
(\% i12)
\end{minipage}
\begin{minipage}[t]{\textwidth}\begin{verbatim}
K0(q1,p1,q2,p2):=(p1^2+q1^2)/2-(p2^2+4*q2^2)/2$\end{verbatim}
\end{minipage}

\noindent
\begin{minipage}[t]{8ex}\bf
(\% i13)
\end{minipage}
\begin{minipage}[t]{\textwidth}\begin{verbatim}
K1(q1,p1,q2,p2):=(q1+q2)^4/4$\end{verbatim}
\end{minipage}

\noindent
\begin{minipage}[t]{8ex}\bf
(\% i14)
\end{minipage}
\begin{minipage}[t]{\textwidth}\begin{verbatim}
K2(q1,p1,q2,p2):=0$\end{verbatim}
\end{minipage}\\

\noindent and then compute the second-order normal form 
(here and below, the Maxima output has been slightly edited in order to fit the page):

\noindent
\begin{minipage}[t]{8ex}\bf
(\% i15)
\end{minipage}
\begin{minipage}[t]{\textwidth}\begin{verbatim}
pnormal2(K0,K1,K2,%lambda);\end{verbatim}
\end{minipage}
\begin{align*}
& \frac{\lambda^2}{4587520} \left(255168 q_2^6+\left( 1184256 q_1^2+191376 p_2^2+1184256 p_1^2\right)\,q_2^4\right. \\
&+ \left. \left( 225792 {q_1^{4}}+\left( 592128 {p_2^{2}}+4580352 {p_1^{2}}\right) \, {q_1^{2}}+47844 {p_2^{4}}+592128 {p_1^{2}}\, {p_2^{2}}+225792 {p_1^{4}}\right) \, {q_2^{2}}\right.\\
&+\left.\left( 2064384 p_1\, p_2\, {q_1^{3}}-2064384 {p_1^{3}}\, p_2\, q_1\right) \, q_2+48384 {q_1^{6}}+\left( 314496 {p_2^{2}}+145152 {p_1^{2}}\right) \, {q_1^{4}}\right.\\
&+\left.\left( 74016 {p_2^{4}}-403200 {p_1^{2}}\, {p_2^{2}}+145152 {p_1^{4}}\right) \, {q_1^{2}}+3987 {p_2^{6}}+74016 {p_1^{2}}\, {p_2^{4}}+314496 {p_1^{4}}\, {p_2^{2}}+48384 {p_1^{6}}\right) \\
&+\frac{\lambda}{512}  \left( 48 {q_2^{4}}+\left( 192 {q_1^{2}}+24 {p_2^{2}}+192 {p_1^{2}}\right) \, {q_2^{2}}+48 {q_1^{4}}+\left( 48 {p_2^{2}}+96 {p_1^{2}}\right) \, {q_1^{2}}+3 {p_2^{4}}+48 {p_1^{2}}\, {p_2^{2}}+48 {p_1^{4}}\right) \\
&-\frac{4 {q_2^{2}}+{p_2^{2}}}{2}+\frac{{q_1^{2}}+{p_1^{2}}}{2}\mbox{}\tag{\% o15} \label{o15}
\end{align*}

\noindent Now, the term $N_2$ can be easily extracted:

\noindent
\begin{minipage}[t]{8ex}\bf
(\% i16)
\end{minipage}
\begin{minipage}[t]{\textwidth}\begin{verbatim}
expand(coeff(%,%lambda,2));\end{verbatim}
\end{minipage}
\begin{align*}
&\frac{3987 {{q_2}^{6}}}{71680}+\frac{2313 {{q_1}^{2}}\, {{q_2}^{4}}}{8960}+\frac{11961 {{p_2}^{2}}\, {{q_2}^{4}}}{286720}+\frac{2313 {{p_1}^{2}}\, {{q_2}^{4}}}{8960}+\frac{63 {{q_1}^{4}}\, {{q_2}^{2}}}{1280}+\frac{2313 {{p_2}^{2}}\, {{q_1}^{2}}\, {{q_2}^{2}}}{17920}\\
&+\frac{639 {{p_1}^{2}}\, {{q_1}^{2}}\, {{q_2}^{2}}}{640}+\frac{11961 {{p_2}^{4}}\, {{q_2}^{2}}}{1146880}+\frac{2313 {{p_1}^{2}}\, {{p_2}^{2}}\, {{q_2}^{2}}}{17920}+\frac{63 {{p_1}^{4}}\, {{q_2}^{2}}}{1280}+\frac{9 p_1\, p_2\, {{q_1}^{3}}\, q_2}{20}\\
&-\frac{9 {{p_1}^{3}}\, p_2\, q_1\, q_2}{20}+\frac{27 {{q_1}^{6}}}{2560}+\frac{351 {{p_2}^{2}}\, {{q_1}^{4}}}{5120}+\frac{81 {{p_1}^{2}}\, {{q_1}^{4}}}{2560}+\frac{2313 {{p_2}^{4}}\, {{q_1}^{2}}}{143360}-\frac{45 {{p_1}^{2}}\, {{p_2}^{2}}\, {{q_1}^{2}}}{512}\\
&+\frac{81 {{p_1}^{4}}\, {{q_1}^{2}}}{2560}+\frac{3987 {{p_2}^{6}}}{4587520}+\frac{2313 {{p_1}^{2}}\, {{p_2}^{4}}}{143360}+\frac{351 {{p_1}^{4}}\, {{p_2}^{2}}}{5120}+\frac{27 {{p_1}^{6}}}{2560}\mbox{}\tag{\% o16} 
\end{align*}

\noindent
\begin{minipage}[t]{8ex}\bf
(\% i17)
\end{minipage}
\begin{minipage}[t]{\textwidth}\begin{verbatim}
define(N2(q1,p1,q2,p2),%)$\end{verbatim}
\end{minipage}\\

\noindent And, finally, the reduction to Hopf variables can be achieved as follows:

\noindent
\begin{minipage}[t]{8ex}\bf
(\% i18)
\end{minipage}
\begin{minipage}[t]{\textwidth}\begin{verbatim}
phopf6res12(expr):=block(
[aux,list_coeff,eq,eqs,W,Wp,U,Up,a,l,
w:[q1^2+p1^2,4*q2^2+p2^2],
u:[-4*p1*q1*q2-p2*q1^2+p1^2*p2,
   -2*q1^2*q2+2*p1^2*q2+2*p1*p2*q1]],
W:makelist(w[1]^i*w[2]^(3-i),i,makelist(j,j,0,3)),
Wp:makelist(%rho[1]^i*%rho[2]^(3-i),i,makelist(j,j,0,3)),
U:makelist(u[1]^i*u[2]^(2-i),i,makelist(j,j,0,2)),
Up:makelist(%rho[3]^i*%rho[4]^(2-i),i,makelist(j,j,0,2)),   
a:makelist(a[k],k,1,length(W)+length(U)),
aux:facsum(expandwrt(expr-sum(a[i]*W[i],i,1,length(W))-
                     sum(a[i+length(W)]*U[i],i,1,length(U)),
                     q1,p1,q2,p2),
           q1,p1,q2,p2),
list_coeff:coeffs(aux,q1,p1,q2,p2),
l:length(list_coeff),
for j:2 thru l do (k:j-1, eq[k]:first(list_coeff[j])),
eqs:makelist(eq[k],k,1,l-1),
subst(first(algsys(eqs,a)),
      sum(a[i]*Wp[i],i,1,length(Wp))
      +sum(a[i+length(W)]*Up[i],i,1,length(U)))
)$\end{verbatim}
\end{minipage}

\noindent
\begin{minipage}[t]{8ex}\bf
(\% i19)
\end{minipage}
\begin{minipage}[t]{\textwidth}\begin{verbatim}
phopf6res12(N2(q1,p1,q2,p2));\end{verbatim}
\end{minipage}
\begin{align*}
\tag{\% o19} 
&-\frac{{{\rho }_{4}^{2}}\, \left( 5120 \% r_1-63\right) }{5120}-\frac{{{\rho }_{3}^{2}}\, \left( 5120 \% r_1-351\right) }{5120}+{{\rho }_{1}^{2}}\, {{\rho }_2}\,\% r_1 \\
&+\frac{3987 {{\rho }_{2}^{3}}}{4587520}+\frac{2313 {{\rho }_1}\, {{\rho }_{2}^{2}}}{143360}+\frac{27 {{\rho }_{1}^{3}}}{2560}\mbox{}
\end{align*}

\noindent
\begin{minipage}[t]{8ex}\bf
(\% i20)
\end{minipage}
\begin{minipage}[t]{\textwidth}\begin{verbatim}
subst(%r1=0,expand(%));\end{verbatim}
\end{minipage}
\[\displaystyle
\tag{\% o20} 
\frac{63 {{\rho }_{4}^{2}}}{5120}+\frac{351 {{\rho }_{3}^{2}}}{5120}+\frac{3987 {{\rho }_{2}^{3}}}{4587520}+\frac{2313 {{\rho }_1}\, {{\rho }_{2}^{2}}}{143360}+\frac{27 {{\rho }_{1}^{3}}}{2560}\mbox{}
\]\\
Of course, a similar analysis can be done for $N_1$. The resulting expression appears in
\cite{AVV17} (see equation (16) in that paper), applied to the determination of
the existence of closed, stable orbits in the Pais-Uhlenbeck oscillator.

\section{Numerical study: Poincar\'e sections}\label{sec3}

Given a Hamiltonian $H\in\mathcal{C}^\infty (\mathbb{R}^{2n})$, the Maxima package
\texttt{poincare} provides several functions to study its Poincar\'e sections.
The functionality of this package, and even the syntax, is similar to the
package \texttt{DEtools} in Maple\texttrademark\, but it offers two advantages: first,
it uses free (both as in `freedom' and as in `free beer') software and, second,
it is almost three times faster, thus being a serious competitor for long
computations. Moreover, as we will see in Section \ref{sec4}, in conjunction with
\KeTCindy\ animation movies describing the evolution of the system in phase space can
be easily constructed.

It should be stressed that Maxima is a CAS, not a language intended for numerical computations such
as Octave. Thus, speed in computations is not one of its goals, nor was it designed to achieve it.
However, the fact that LISP is its underlying programming language, allows the possibility of
writing specialized routines using declared variables, that can be then compiled. The package
\texttt{poincare} uses a compiled version of the Runge-Kutta method, called \texttt{rkfun},
developed by Richard Fateman (\url{http://people.eecs.berkeley.edu/~fateman/lisp/rkfun.lisp}),
and this is the ultimate reason for the gain in speed.

The function called \texttt{hameqs} constructs the Hamiltonian equations for a given Hamiltonian $H(q1,p1,\ldots ,q2,p2)$. Any
names can be used for the variables, but they must be given in pairs
``coordinate, conjugate momentum''. A good choice (used internally) is 
$(q1,p1,...,qn,pn)$. A name must be
provided for the components of the Hamiltonian vector field
$$
X_H(q_1,p_1,...,q_n,p_n)=\sum^n_{i=1}\left( 
\frac{\partial H}{\partial p_i}\frac{\partial}{\partial q_i}
-\frac{\partial H}{\partial q_i}\frac{\partial}{\partial p_i}
\right)\,.
$$
Once a name, say $XH$, is chosen, the components of the Hamiltonian vector field
will be globally defined functions $XHj$ with $1\leq j\leq 2n$, where $n$ is the
number of degrees of freedom, and will be available to Maxima. Notice that, for
instance,
$$
XH1(t,q_1,p_1,...,q_n,p_n)=\frac{\partial H}{\partial p_1}
$$
and
$$
XH2(t,q_1,p_1,...,q_n,p_n)=-\frac{\partial H}{\partial q_1}\,.
$$
Although we will work with \emph{autonomous}
Hamiltonian systems, the components $XHj$ returned by this command will
have the set $(t,q_1,p_1,...,q_n,p_n)$ as arguments. This is necessary 
to maintain consistency with the \texttt{rkfun} routine (implementing the
Runge-Kutta method), which can work
with both, autonomous and non-autonomous systems.

The function \texttt{poincare3d} constructs the projection of the Hamiltonian 
orbits along a certain coordinate which is given as an argument \texttt{coord}.
Other arguments are: a list
of initial conditions $\mathtt{inicond} =[q_1(0),p_1(0),\ldots ,q_n(0),p_n(0)]$, and
a list characterizing the time domain $\mathtt{timestep}=[t,t_{ini},t_{fin},step]$.
Thus, the syntax is \texttt{poincare3d(H,name,inicond,timestep,coord)}. The package is loaded
with

\noindent
\begin{minipage}[t]{8ex}\bf
(\%{}2i1) 
\end{minipage}
\begin{minipage}[t]{\textwidth}\begin{verbatim}
batch("poincare.mac")$
\end{verbatim}
\end{minipage}

As a simple example, let us construct the $3D-$surface of a couple of harmonic oscillators
(this is based on Chapter 9 of \cite{Lyn01}, where a similar discussion using Maple\texttrademark\ is
presented):

\noindent
\begin{minipage}[t]{8ex}\bf
(\%{}i22) 
\end{minipage}
\begin{minipage}[t]{\textwidth}\begin{verbatim}
H(x,v,q,p):=w1*(x^2+v^2)/2+w2*(q^2+p^2)/2$
\end{verbatim}
\end{minipage}

\noindent The corresponding Hamiltonian equations are:

\noindent
\begin{minipage}[t]{8ex}\bf
(\%{}i23) 
\end{minipage}
\begin{minipage}[t]{\textwidth}\begin{verbatim}
hameqs(H,XH);
\end{verbatim}
\end{minipage}
\[\displaystyle
\tag{\%{}o7}\label{o23} 
[v\,\mathit{w1},-\mathit{w1}x,p\,\mathit{w2},-q\,\mathit{w2}]\mbox{}
\]

\noindent and we fix the values of the frequencies $w1=1$, $w2=3$ so they are commensurable:

\noindent
\begin{minipage}[t]{8ex}\bf
(\%{}24) 
\end{minipage}
\begin{minipage}[t]{\textwidth}\begin{verbatim}
[w1,w2]:[8,3]$
\end{verbatim}
\end{minipage}

Next we plot the $3D$ Poincar\'e surface by projecting along the $p$ coordinate
(thus, the resulting graphics has $(x,v,q)$ coordinates):

\noindent
\begin{minipage}[t]{8ex}\bf
(\%{}i25) 
\end{minipage}
\begin{minipage}[t]{\textwidth}\begin{verbatim}
data1:poincare3d(H,XH,[0.3,0.5,0,1.5],[t,0,40,0.01],p)$
\end{verbatim}
\end{minipage}

\noindent
\begin{minipage}[t]{8ex}\bf
(\%{}i26) 
\end{minipage}
\begin{minipage}[t]{\textwidth}\begin{verbatim}
draw3d(title="Poincare section in 3D",
dimensions=[350,500],view=[85,30],
xlabel="x",ylabel="v",zlabel="q",
xtics=1,ytics=1,
surface_hide=true,color="light-blue",
    explicit(0,x,-1.35,1.35,y,-1.35,1.35),
point_size=0,points_joined=true,color=black,line_width=1,
    points(data1),
user_preamble="set xyplane at -1.8",color="light-blue",
    explicit(-1.8,u,-1.5,1.35,v,-1.35,1.35),
point_size=1,point_type=filled_circle,color=red,points_joined=false,
    points([[-0.56,0,-1.78],[0.28,-0.52,-1.78],[0.3,0.48,-1.78],
            [-0.56,0,0],[0.28,-0.52,0],[0.3,0.48,0]]));
\end{verbatim}
\end{minipage}

\[\displaystyle
\tag{\%{}t26,\%{}o26}\label{t26} 
\includegraphics[scale=0.37]{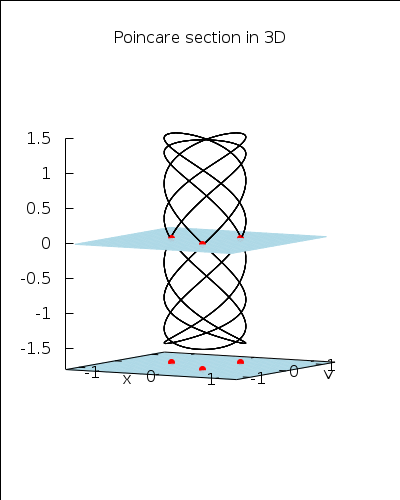}
\]

The function \texttt{poincare2d} constructs the surface of section selected by a list
of arguments of the form $\mathtt{scene}=[q0,c,qi,qj]$, that is, the surface 
$q0=c$  in which coordinates  $[qi,qj]$ are shown.
The method used in the computation of the Poincar\'e surface is that described in
the paper \cite{Cheb96} we select a set of initial conditions,
follow the corresponding orbit numerically, and detect where we have crossed the $q0=c$ surface by
looking at changes of sign in the list of values for this coordinate minus $c$.
In the previous example, we plotted the $3D-$Poincar\'e surface of a couple of
commensurable oscillators, and we included a $2D-$section (corresponding to $q=0$)
showing that the periodicity
of the system reflects itself in the discrete character of the $2D-$Poincar\'e map
(only three points appear in it). Now we can check this directly with \texttt{poincare2d} (notice the selection of the $q=0$ section in the last argument,
\texttt{[q,0,x,v]}):

\noindent
\begin{minipage}[t]{8ex}\bf
(\%{}i27) 
\end{minipage}
\begin{minipage}[t]{\textwidth}\begin{verbatim}
data2:poincare2d(H,XH,[0.3,0.5,0,1.5],[t,0,40,0.01],[q,0,x,v])$
\end{verbatim}
\end{minipage}

\noindent
\begin{minipage}[t]{8ex}\bf
(\%{}i28) 
\end{minipage}
\begin{minipage}[t]{\textwidth}\begin{verbatim}
draw2d(title="Poincare section in 2D",
xlabel="x",ylabel="v",
xtics=0.2,
point_size=1,point_type=7,color=red,
points_joined=false,proportional_axes=xy,
    points(data2));
\end{verbatim}
\end{minipage}

\[\displaystyle
\tag{\%{}t28,\%{}o28}\label{t28} 
\includegraphics[scale=0.37]{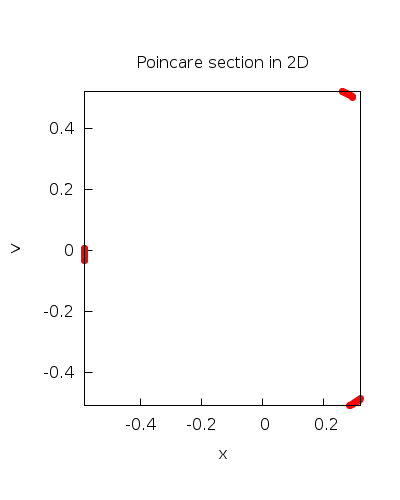}\mbox{}
\]


For a different example, let us consider the case of a elastic
pendulum \cite{Car94}, with Hamiltonian

\noindent
\begin{minipage}[t]{8ex}\bf
(\%{}i29) 
\end{minipage}
\begin{minipage}[t]{\textwidth}\begin{verbatim}
H(q1,p1,q2,p2):=(p1^2+p2^2)/2+(q1^2+q2^2)/2-0.75*q1^2*(1+q2)/2$
\end{verbatim}
\end{minipage}

We can obtain an analytic expression for $p_2$ once the energy $E$, and the initial values of $(q_1,p_1,q_2)$ are known. Here we work with $E=0{.}00875$:

\noindent
\begin{minipage}[t]{8ex}\bf
(\%{}i30) 
\end{minipage}
\begin{minipage}[t]{\textwidth}\begin{verbatim}
solve(H(q1,p1,q2,p2)=0.00875,p2);
\end{verbatim}
\end{minipage}
 
\begin{align*}
[p_2=-\frac{\sqrt{-4{{q_2}^{2}}+3{{q_1}^{2}}\,
q_2-{{q_1}^{2}}-4{{p_1}^{2}}+0.07}}{2},
p_2=\frac{\sqrt{-4{{q_2}^{2}}+3{{q_1}^{2}}\,
q_2-{{q_1}^{2}}-4{{p_1}^{2}}+0.07}}{2}]
\end{align*}\[\displaystyle
\tag{\%{}o30}\label{o30}\]

Let us define the corresponding functions:

\noindent
\begin{minipage}[t]{8ex}\bf
(\%{}i31) 
\end{minipage}
\begin{minipage}[t]{\textwidth}\begin{verbatim}
define(f(q1,p1,q2),rhs(first(%)))$
\end{verbatim}
\end{minipage}

\noindent
\begin{minipage}[t]{8ex}\bf
(\%{}i32) 
\end{minipage}
\begin{minipage}[t]{\textwidth}\begin{verbatim}
define(g(q1,p1,q2),rhs(second(%th(2))))$
\end{verbatim}
\end{minipage}

Now, we compute the $q_2=0$ surface of section, for enough initial conditions $(q_1,p_1,q_2)$ ($10$ different sets)  using the positive value of $p_2$, and joining all the resulting points in a
big list of $2D-$coordinates called \texttt{points1}:

\noindent
\begin{minipage}[t]{8ex}\bf
(\%{}i33) 
\end{minipage}
\begin{minipage}[t]{\textwidth}\begin{verbatim}
for j:1 thru 10 do data1[j]:poincare2d(H,XH,
[0.15,j/100,0.001,g(0.15,j/100,0.001)],[t,0,1000,0.01],[q2,0,q1,p1])$
\end{verbatim}
\end{minipage}

For future reference, here is the time invested in the computation:

\noindent
\begin{minipage}[t]{8ex}\bf
(\%{}i34) 
\end{minipage}
\begin{minipage}[t]{\textwidth}\begin{verbatim}
time(%);
\end{verbatim}
\end{minipage}
\[\displaystyle
\tag{\%{}o61}\label{o61} 
[17.222]\mbox{}
\]

\noindent
\begin{minipage}[t]{8ex}\bf
(\%{}i35) 
\end{minipage}
\begin{minipage}[t]{\textwidth}\begin{verbatim}
points1:xreduce(append,create_list(data1[j],j,makelist(k,k,1,10)))$
\end{verbatim}
\end{minipage}

The following figure is the plot of these points on the Poincar\'e surface:

\noindent
\begin{minipage}[t]{8ex}\bf
(\%{}i36) 
\end{minipage}
\begin{minipage}[t]{\textwidth}\begin{verbatim}
draw2d(title="Poincare sections E=0.00875",
xlabel="q1",ylabel="p1",
point_type=7,point_size=0.1,
    points(points1)
);
\end{verbatim}
\end{minipage}
\[\displaystyle
\tag{\%{}t36,\%{}o36}\label{t36} 
\includegraphics[width=.67\linewidth,height=.52\textheight,keepaspectratio]{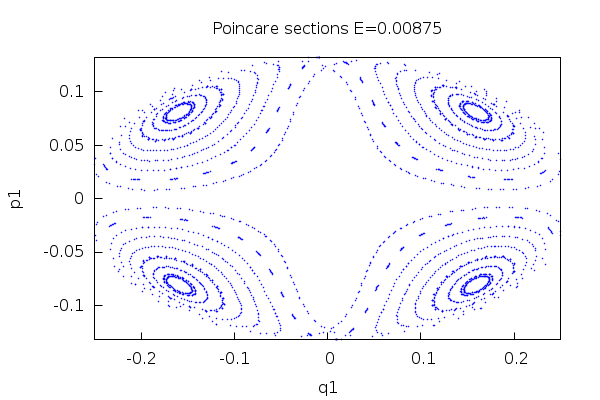}\mbox{}\]

In order to complete the section, we must select another set of initial conditions, whose orbits pass through the empty region at the center.
This time we use negative values of the momentum $p_2$:

\noindent
\begin{minipage}[t]{8ex}\bf
(\%{}i37) 
\end{minipage}
\begin{minipage}[t]{\textwidth}\begin{verbatim}
for j:1 thru 10 do data2[j]:poincare2d(H,XH,
[2*j/100,0,j/100+0.0025,f(2*j/100,0,j/100+0.0025)],[t,0,1000,0.01],
[q2,0,q1,p1])$
\end{verbatim}
\end{minipage}

\noindent
\begin{minipage}[t]{8ex}\bf
(\%{}i38) 
\end{minipage}
\begin{minipage}[t]{\textwidth}\begin{verbatim}
time(%);
\end{verbatim}
\end{minipage}
\[\displaystyle
\tag{\%{}o38}\label{o38} 
[17.369]\mbox{}
\]

The $2D-$dimensional coordinates of the corresponding points are stored in the list \texttt{points2}
and then plotted with the aid of the \texttt{draw2d} command, which admits lots of optional
arguments to fine tuning the appearance of the figure. Here, we specify that points be represented by
filled circles (\texttt{point\_type=7}) with a given radius size (\texttt{point\_size=0.1}): 

\noindent
\begin{minipage}[t]{8ex}\bf
(\%{}i39) 
\end{minipage}
\begin{minipage}[t]{\textwidth}\begin{verbatim}
points2:xreduce(append,create_list(data2[j],j,makelist(k,k,1,10)))$
\end{verbatim}
\end{minipage}

\noindent
\begin{minipage}[t]{8ex}\bf
(\%{}i40) 
\end{minipage}
\begin{minipage}[t]{\textwidth}\begin{verbatim}
draw2d(title="Poincare sections E=0.00875",
xlabel="q1",ylabel="p1",
point_type=7,point_size=0.1,
    points(points2)
);
\end{verbatim}
\end{minipage}
\[\displaystyle
\tag{\%{}t40,\%{}o40}\label{t40} 
\includegraphics[width=.67\linewidth,height=.52\textheight,keepaspectratio]{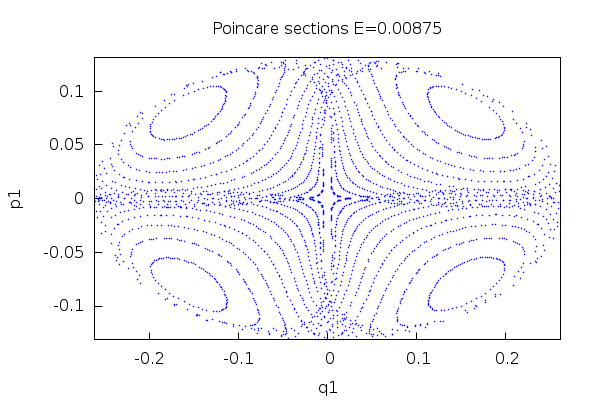}\mbox{}\]

The full Poincar\'e section is obtained by joining both sets of points. The Maxima command
\texttt{append} does exactly that when two lists are given:

\noindent
\begin{minipage}[t]{8ex}\bf
(\%{}i41) 
\end{minipage}
\begin{minipage}[t]{\textwidth}\begin{verbatim}
draw2d(title="Poincare sections E=0.00875",
xlabel="q1",ylabel="p1",
point_type=7,point_size=0.1,
    points(append(points1,points2))
);
\end{verbatim}
\end{minipage}
\[\displaystyle
\tag{\%{}t41}\label{t41} 
\includegraphics[width=.75\linewidth,height=.60\textheight,keepaspectratio]{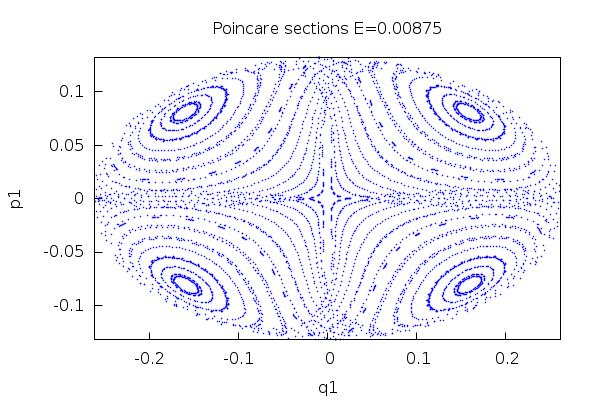}\mbox{}\]

We have done our computations with a fixed value for the energy
$E=0{.}00875$. The same steps can be followed to consider other
values. The collection of Poincar\'e sections so obtained is a
valuable tool to visualize the dynamics of the system.
In Section \ref{sec4} we will see how to put all these
sections together in the form of a movie animating the evolution
in phase space.

\section{Graphical study with KeTCindy}\label{sec4}

\ketpic{} is a macro package involving several mathematical software
for producing high-quality figures to be inserted into \LaTeX\ documents, developed by
one of the authors (ST). It can use the DGS Cinderella
as a graphical interface through \ketcindy{}, another set of macros which acts as a interface
between them.
The reason for choosing Cinderella is that it has its own scripting language, CindyScript,
featuring a simple to understand syntax.
Using CindyScript, we have added a layer to \ketcindy{} to call other software such as Maxima,
\textbf{\textsf{R}}, Fricas, Risa/Asir or C. We refer to previous works for more details on
the CindyScript syntax \cite{icms2016,cas2016,iccsa2017,castr2017}; here we proceed in a more
direct way, explaining how \ketcindy\ calls to Maxima 
by way of an example devoted to find the indefinite and definite integrals of a function.
For this, the following code must be inserted into the script editor of Cinderella:

\begin{verbatim}
cmdL=[
 "f(x):=sin(x)+cos(2*x)",[],
 "ans1:integrate",["f(x)","x"],
 "ans2:integrate",["f(x)","x",0,"%pi/3"],
 "ans1::ans2",[]
];
CalcbyM("ans",cmdL,[""]);
\end{verbatim}

Here \texttt{cmdL} is a list of Maxima commands which are parsed sequentially. By executing 
\texttt{CalcbyM} in the next line, a file called \texttt{simpleexampleans.txt}
in the directory \texttt{ketwork} with the contents given below will be created:

\begin{verbatim}
writefile("ketwork/simpleexampleans.txt")$/*##*/
powerdisp:false$/*##*/
display2d:false$/*##*/
linel:1000$/*##*/
f(x):=sin(x)+cos(2*x)$/*##*/
ans1:integrate(f(x),x)$/*##*/
ans2:integrate(f(x),x,0,%pi/3)$/*##*/
disp(ans1)$/*##*/
disp(ans2)$/*##*/
closefile()$/*##*/
quit()$/*##*/
\end{verbatim}

These instructions will be processed by Maxima, and the results will be passed to
\ketcindy{} to generate either direct graphical output in Cinderella or a \LaTeX\ file with
the corresponding code to generate the graphics that
can be inserted in another document. In more detail, \texttt{calcbyM} will sequentially do the following:

\begin{enumerate}
\item Create the \texttt{txt} file to be processed by Maxima.
\item Create a batch file \texttt{kc.sh(bat)} to call Maxima.
\item Call a Java program to execute the batch file above.
\item Hand the result from Maxima, parsed as strings, to \ketcindy{}.
\end{enumerate}

In the above example, the result \texttt{ans} is a list containing two strings:
\begin{verbatim}
[sin(2*x)/2-cos(x),(sqrt(3)+2)/4] 
\end{verbatim}

This result can be directly used in \ketcindy; for example, we can draw the graph of the indefinite integral with the following command:
\begin{verbatim}
Plotdata(``1'',ans_1,''x'');
\end{verbatim}

\begin{center}
{\unitlength=1cm%
\begin{picture}%
(10,4)(-5,-2)%
\linethickness{0.008in}
\polyline(-5.00000,-0.01165)(-4.80000,-0.00034)(-4.60000,0.00071)(-4.40000,0.01487)%
(-4.20000,0.06296)(-4.00000,0.15896)(-3.80000,0.30701)(-3.60000,0.49992)(-3.40000,0.71974)%
(-3.20000,0.94002)(-3.00000,1.12970)(-2.80000,1.25786)(-2.60000,1.29862)(-2.40000,1.23548)%
(-2.20000,1.06430)(-2.00000,0.79455)(-1.80000,0.44846)(-1.60000,0.05839)(-1.40000,-0.33746)%
(-1.20000,-0.70009)(-1.00000,-0.99495)(-0.80000,-1.19649)(-0.60000,-1.29136)(-0.40000,-1.27974)%
(-0.20000,-1.17478)(0.00000,-1.00000)(0.20000,-0.78536)(0.40000,-0.56238)(0.60000,-0.35932)%
(0.80000,-0.19692)(1.00000,-0.08565)(1.20000,-0.02463)(1.40000,-0.00247)(1.60000,0.00001)%
(1.80000,0.00594)(2.00000,0.03775)(2.20000,0.11270)(2.40000,0.23931)(2.60000,0.41516)%
(2.80000,0.62659)(3.00000,0.85028)(3.20000,1.05657)(3.40000,1.21385)(3.60000,1.29359)%
(3.80000,1.27493)(4.00000,1.14832)(4.20000,0.91756)(4.40000,0.59979)(4.60000,0.22360)%
(4.80000,-0.17466)(5.00000,-0.55567)%
\polyline(-5.00000,0.00000)(5.00000,0.00000)%
\polyline(0.00000,-2.00000)(0.00000,2.00000)%
\settowidth{\Width}{$x$}\setlength{\Width}{0\Width}%
\settoheight{\Height}{$x$}\settodepth{\Depth}{$x$}\setlength{\Height}{-0.5\Height}\setlength{\Depth}{0.5\Depth}\addtolength{\Height}{\Depth}%
\put(5.0500000,0.0000000){\hspace*{\Width}\raisebox{\Height}{$x$}}%
\settowidth{\Width}{$y$}\setlength{\Width}{-0.5\Width}%
\settoheight{\Height}{$y$}\settodepth{\Depth}{$y$}\setlength{\Height}{\Depth}%
\put(0.0000000,2.0500000){\hspace*{\Width}\raisebox{\Height}{$y$}}%
\settowidth{\Width}{O}\setlength{\Width}{-1\Width}%
\settoheight{\Height}{O}\settodepth{\Depth}{O}\setlength{\Height}{-\Height}%
\put(-0.0500000,-0.0500000){\hspace*{\Width}\raisebox{\Height}{O}}%
\end{picture}}%
\end{center}

As mentioned, \ketcindy\ can also produce a \TeX\ animation. We illustrate this feature
with the Poincar\'e sections of an elastic pendulum as an example. We begin by defining
\texttt{Elist} as a list of increasing energies:
\begin{verbatim}
Elist=[0.00875,0.0125,0.01625,0.02,0.02375,0.0275,0.03125,0.035,0.03875];
\end{verbatim}

Next, we generate the corresponding data for each energy using Maxima with the package
\texttt{poincare}. The following list of Maxima commands is just the same
as the one described in Section \ref{sec3}:

\begin{verbatim}
cmdL1=concat(Mxload("rkfun.lisp"),Mxbatch("pdynamics.mac"));
cmdL1=concat(cmdL1,Mxbatch("poincare.mac"));
cmdL1=concat(cmdL1,[
 "H(q1,p1,q2,p2):=(p1^2+p2^2)/2+(q1^2+q2^2)/2-0.75*q1^2*(1+q2)/2",[],
 "ans:solve(H(q1,p1,q2,p2)=E,p2)",[],
 "define(f(q1,p1,q2),rhs(first(%)))",[],
 "define(g(q1,p1,q2),rhs(second(%th(2))))",[]
]);
forall(1..9,nn,
  cmdL2=[
   "E:"+textformat(Elist_nn,6),[],
   "for j:1 thru 10 do data1[j]:poincare2d(H,XH,[0.15,j/100,0.001,
        g(0.15,j/100,0.001)],[t,0,1000,0.01],[q2,0,q1,p1])",[],
   "points1:xreduce(append,create_list(data1[j],j,makelist(k,k,1,10)))",[],
   "for j:1 thru 10 do data2[j]:poincare2d(H,XH,[2*j/100,0,j/100+0.0025,
        f(2*j/100,0,j/  100+0.0025)],[t,0,1000,0.01],[q2,0,q1,p1])",[],
   "points2:xreduce(append,create_list(data2[j],j,makelist(k,k,1,10)))",[],
   "points1::points2",[]
  ];
  cmdL=concat(cmdL1,cmdL2);
  CalcbyM("Points",cmdL,[mr,"Wait=40"]);
);
\end{verbatim}

Each result is stored in a text file (with extension \texttt{txt}) with its name constructed appending 
a number sequentially to the prefix \texttt{ptdata}.
Finally, we define a function \texttt{mf(nn)}, which describes the animation frame numbered 
\texttt{nn}.

\begin{verbatim}
mf(nn):=(
 regional(tmp,Points1,Points2);
 Com1st("ReadOutData('ptdata"+text(nn)+".txt')");
 Setcolor("red");
 Pointdata("1","Points",[red,"Size=0.4"]);
 Setcolor("black");
 Expr(D,"c","E="+textformat(Elist_(nn),6));
);
Setpara("poincare","mf(nn)",1..9,
   ["m","Frate=3","Scale=0.7","OpA=[loop]"]);
\end{verbatim}

The animation is generated by putting together all the frames with the command \texttt{Mkanimation()}.
The result, shown below, requires Adobe Acrobat Reader\texttrademark\ for playback\footnote{Currently, it is the only PDF reader capable of that. Some versions of the KDE reader Okular have been reported to be able of
reproducing some animations, but we have not had success when using it.}: 

\begin{center}
\begin{animateinline}[autoplay,loop,controls]{3}%
\scalebox{0.7}{\input{p001.tex}}%
\newframe%
\scalebox{0.7}{\input{p002.tex}}%
\newframe%
\scalebox{0.7}{\input{p003.tex}}%
\newframe%
\scalebox{0.7}{\input{p004.tex}}%
\newframe%
\scalebox{0.7}{\input{p005.tex}}%
\newframe%
\scalebox{0.7}{\input{p006.tex}}%
\newframe%
\scalebox{0.7}{\input{p007.tex}}%
\newframe%
\scalebox{0.7}{\input{p008.tex}}%
\newframe%
\scalebox{0.7}{\input{p009.tex}}%
\end{animateinline}
\end{center}
\medskip

This graphical analysis illustrates several features common to any perturbed system of the form
$H=H_0+\varepsilon H_1$, where $H_0$ is an integrable Hamiltonian and $\varepsilon \sim 0$. Starting
at low values of the energy, many closed curves can be detected, corresponding to periodic motions
of the system. Those closed curves are intersections of the tori determined by the non-perturbed
part with the Poincar\'e surface. As the energy of the systems increases, the tori are destroyed and 
the trajectories initially confined to them start wandering all over the phase space. At a certain
point, we can not distinguish any periodicity and the behavior is completely chaotic. This generic
picture is the content of the famous KAM (for Kolgomorov, Arnold and Moser) theorem, although this
theorem refers to increasing values of the perturbation parameter rather than the total energy (see
\cite{Cue92} for the application to this case, along with some comments on the applicability of the 
KAM theorem, which is not immediate).  

\section{Conclusions}

The symbolic computation of second-order normal form for perturbed Hamiltonian
systems can be quickly computed in closed form with the aid of the Maxima CAS,
directly in terms of the Hopf invariants. The package \texttt{pdynamics} shows
a practical implementation.

The Maxima package \texttt{poincare} can reproduce the results appearing in
textbooks and research papers dealing with Hamiltonian systems.
The graphical output quality is quite good, comparable (to say the least) to that
of commercial software, but at no cost (for comparison, Maple\texttrademark\, in
its student's version costs $1\,000$USD.) Regarding computation times, the Maxima
version outperforms commercial competitors: the heaviest computation
in this paper is executed in \eqref{o38} while, for instance, 
the same takes $50$ seconds in Maple\texttrademark\ \footnote{Used here: 
Maple 2016:1a (build 1133417).  Maplesoft, a division of Waterloo Maple Inc., 
Waterloo, Ontario. Maxima version was 5.38.0.} (as can be seen in the worksheet
\url{http://galia.fc.uaslp.mx/~jvallejo/ElasticPendulum-MapleSession.pdf}, for
which the same computer was used). Maxima only requires a third of this time.

On the other hand, \ketcindy\ in combination with Maxima can produce \TeX\ animations,
ready for use in complex documents which require high-quality graphics, such as
research papers or handouts to be used in teaching.
The union of these features results in an easy-to-use,
powerful integrated system particularly suitable for studying the dynamics of
Hamiltonian systems.


\subsection*{Acknowledgment}
The authors express their gratitude to Masataka Kaneko (T\={o}h\={o} University) and 
Yasuyuki Nakamura (Nagoya Institute of Technology)
 for many fruitful discussions about the topics in this paper. Thanks are also due to
 Richard Fateman for developing \texttt{rkfun}.
\end{document}